# New orthogonal hybrid function based numerical method to solve system of fractional order differential equations

Seshu Kumar Damarla[1], Madhusree Kundu[2]

[1]Department of Chemical Engineering, C.V. Raman College of Engineering, Bhubaneswar, Odisha, India.

[1]E-mail: seshukumar.damarla@gmail.com

[2]Department of Chemical Engineering, National Institute of Technology Rourkela, India-769008.

[2]E-mail: mkundu@nitrkl.ac.in

**Abstract**

In this paper, an easy-to-implement and computationally effective numerical method based on the new orthogonal hybrid functions is developed to solve system of fractional order differential equations numerically. The new orthogonal hybrid functions are hybrid of the piecewise constant orthogonal sample-and-hold functions and the piecewise linear orthogonal right-handed triangular functions. The proposed method uses the generalized one-shot operational matrices which approximate the Riemann-Liouville fractional order integral in the orthogonal hybrid function domain. The convergence of the numerical method is studied. Illustrative examples such as fractional order smoking model, fractional order model for lung cancer, fractional order model of Hepatitis B infection etc. are solved by the proposed numerical method. The results prove the validity and reliability of the proposed numerical method.

## 1. Introduction

Fractional calculus is a three centuries old topic like its counterpart integer order calculus but unlike its counterpart it is up till now a fresh subject. It is one of attractive and rapidly growing fields of mathematics which deals with the theory of differentiation and integration of arbitrary order (order can be integer, non-integer or complex number). It has been emerging as an unquestionably crucial subject for applied mathematicians and applied scientists for better comprehending numerous physical processes in varied applied areas of science and engineering such as electrochemistry, physics, geology, astrophysics, seismic wave analysis, sound wave propagation, psychology and life sciences, biology, etc. (Oldham and Spanier (1974), Atanackovic (2014), Hermann (2011), Losa et al. (2005), Mainardi (2010), Golmankhneh (2015), Biyajima et al. (2015), Ahmad and El-Khazali (2007), Song et al. (2010)). Gaining insights into the mentioned physical processes demands the exact solutions of fractional order differential equations representing them. Nonetheless it is extremely difficult (even impossible in some circumstances) to find analytical solution of those fractional order differential equations. This drawback limits the practical applications of the subject. To enable the subject to be successfully applied in the said application areas, some analytical methods (Podlubny (1999)) were proposed but their range of applicability is very narrow so several numerical methods (Diethelm et al. (2002), Odibat and Momani (2008), Adomian (1994), He (1999a), (1999b), Arigoklu and Ozkol (2007), Zurigat et al. (2010)), each having its own advantages and disadvantages, were devised. As is well known that no analytical or numerical method can solve all classes of fractional order differential equations and some real world processes are inherently very complex and highly nonlinear, especially biological processes, that pose computational challenges which the available semi-analytical or numerical methods may not meet. Significant efforts continue to develop more efficient and reliable numerical methods which are able to precisely approximate the solution of fractional order differential equations. From this standpoint, we propose a numerical method to solve system of fractional order ordinary differential equations of the following form.

$${}_{0}^{C}D_{t}^{\alpha}y_{i}(t) = f_{i}(t, y_{1}(t), y_{2}(t), y_{3}(t), \ldots, y_{n}(t)), \ i \in [1, n], \ y_{i}(0) = c_{i}, \ t \in [0,1], \ \alpha \in (0,1]. \qquad (1)$$

Hhere $y_i(t)$ is the $i^{th}$ unknown function, $f_i(t, y_1(t), y_2(t), y_3(t), \ldots, y_n(t))$ can be linear or nonlinear and the fractional derivative, ${}_{0}^{C}D_{t}^{\alpha}$, is the Caputo fractional order derivative

$${}_{0}^{C}D_{t}^{\alpha}f(t) = J^{m-\alpha}D^{m}f(t) = \frac{1}{\Gamma(m-\alpha)}\int_{0}^{t}(t-\tau)^{m-\alpha-1}f^{(m)}(\tau)d\tau, \qquad (2)$$

where $J^{\alpha}f(t)$ is Riemann-Liouville fractional order integral, $J^{\alpha}f(t) = \frac{1}{\Gamma(\alpha)}\int_{0}^{t}(t-\tau)^{\alpha-1}f(\tau)d\tau$.

Our numerical technique is entirely based on the new orthogonal hybrid functions (HF) formed by combining the piecewise constant orthogonal sample-and-hold functions and the piecewise linear orthogonal right-handed triangular functions (Deb et al. (2016)). The new orthogonal hybrid functions were fruitfully employed for analysis and identification of time-invariant, time-varying delay and delay-free systems (Deb et al. (2016)) and network analysis (Sarkar et al. (2014)). In (Deb et al. (2016), (2012)), linear ordinary differential equations of order up to three were only solved using orthogonal HFs but we anticipate that HFs are capable enough to be utilized to solve nonlinear ordinary differential equations, integral equations, integro-differential equations, partial differential equations, stochastic differential equations, etc. Since the integer order calculus is a special case of the fractional calculus and orthogonal HFs possess possible applications in the classical calculus i.e. integer order calculus, we trust that orthogonal HFs may find useful applications even in fractional calculus. To the best of our knowledge, it is the first attempt to use orthogonal HFs to solve complex system of nonlinear fractional order differential equations. The remaining part of the paper is structured as follows. Section 2 introduces orthogonal HFs. The generalized one-shot operational matrices, which are the basis of our numerical method, approximating the Riemann-Liouville fractional order integral in the orthogonal HF domain are derived in Section 3. The numerical method based on the generalized one-shot operational matrices is developed in Section 4. Sections 4 proves theoretically that the proposed numerical method converges the approximate solution of the system of nonlinear fractional order differential equations to its actual solution. The proposed numerical method is tested in Section 6. Section 7 concludes the paper.

## 2. Hybrid functions (HF)

Let us define the $i^{th}$ component of a set of orthogonal hybrid functions containing $m$ component functions in the semi-open interval $[0,T)$ as (Deb et al. (2016))

$$H_i(t) = \begin{cases} a_i S_i(t) + b_i T_i(t), & \text{if } ih \leq t < (i+1)h, \\ 0, & \text{otherwise,} \end{cases}, \quad i = 0,1,2,\ldots,m-1, \tag{3}$$

where $a_i$ and $b_i$ are constants defined as $a_i = c_i$ and $b_i = (c_{i+1} - c_i)$, $S_i(t)$ is the $i^{th}$ sample-and-hold function (SHF) defined as (Deb et al. (2007))

$$S_i(t) = \begin{cases} 1, & \text{if } ih \leq t < (i+1)h \\ 0, & \text{otherwise} \end{cases} \tag{4}$$

and $T_i(t)$ is the $i^{th}$ right-handed triangular function (TF) defined as (Deb et al. (2007))

$$T_i(t) = \begin{cases} \dfrac{t - ih}{h}, & \text{if } ih \leq t < h \\ 0, & \text{otherwise} \end{cases}. \tag{5}$$

Consider a time function, $f(t)$, of Lebesgue measure defined on the interval $[0,T)$. Let us take $m$ equidistant nodes in the given interval with the constant step size $h$ as $t_i = ih$, $i = 0,1,2,\ldots,m$. The function $f(t)$ can be expressed in terms of orthogonal HFs as follows.

$$\begin{aligned} f(t) &\approx \sum_{i=0}^{m-1} H_i(t) = \sum_{i=0}^{m-1} a_i S_i(t) + b_i T_i(t), \\ &= [c_0 \quad c_1 \quad c_2 \quad \cdots \quad \cdots \quad c_{m-1}] S_{(m)}(t) + [(c_1 - c_0) \quad (c_2 - c_1) \quad (c_3 - c_2) \quad \cdots \quad \cdots \quad (c_m - c_{m-1})] T_{(m)}(t), \\ &= C_S^T S_{(m)}(t) + C_T^T T_{(m)}(t), \end{aligned} \tag{6}$$

where

$$C_S^T = [c_0 \quad c_1 \quad c_2 \quad \cdots \quad \cdots \quad c_{m-1}], \quad C_T^T = [(c_1 - c_0) \quad (c_2 - c_1) \quad \cdots \quad \cdots \quad (c_m - c_{m-1})], \quad c_i = f(ih),$$
$$S_{(m)}(t) = [S_0(t) \quad S_1(t) \quad S_2(t) \quad \cdots \quad S_{m-1}(t)], \quad T_{(m)}(t) = [T_0(t) \quad T_1(t) \quad T_2(t) \quad \cdots \quad T_{m-1}(t)].$$

The first order integral of $f(t)$ can be approximated on $[0,T)$ in the orthogonal HF domain as [20]

$$\int_0^t f(\tau)d\tau \approx \int_0^t \left(C_S^T S_{(m)}(\tau) + C_T^T T_{(m)}(\tau)\right)d\tau = C_S^T \int_0^t S_{(m)}(\tau)d\tau + C_T^T \int_0^t T_{(m)}(\tau)d\tau,$$

$$= C_S^T \left(P_{1ss(m)} S_{(m)}(t) + P_{1st(m)} T_{(m)}(t)\right) + C_T^T \left(P_{1ts(m)} S_{(m)}(t) + P_{1tt(m)} T_{(m)}(t)\right), \tag{7}$$

$$= \left(C_S^T P_{1ss(m)} + C_T^T P_{1ts(m)}\right) S_{(m)}(t) + \left(C_S^T P_{1st(m)} + C_T^T P_{1tt(m)}\right) T_{(m)}(t),$$

where

$$P_{1ss(m)} = h[\![0 \ \ 1 \ \ 1 \ \ \cdots \ \ 1]\!], \qquad P_{1st(m)} = h[\![1 \ \ 0 \ \ 0 \ \ \cdots \ \ 0]\!], \qquad P_{1ts(m)} = \frac{h}{2}[\![0 \ \ 1 \ \ 1 \ \ \cdots \ \ 1]\!],$$

$$P_{1tt(m)} = \frac{h}{2}[\![1 \ \ 0 \ \ 0 \ \ \cdots \ \ 0]\!], \ [\![a \ \ b \ \ c]\!] = \begin{bmatrix} a & b & c \\ 0 & a & b \\ 0 & 0 & a \end{bmatrix}.$$

The complementary operational matrices, $P_{1ss(m)}$, $P_{1st(m)}$, $P_{1ts(m)}$ and $P_{1tt(m)}$, are acting as a first order integrator in the orthogonal HF domain. The orthogonal and operational properties of HFs are given Appendix 1.

## 3. Generalized one-shot operational matrices for fractional order integral of $f(t)$

In this section, we shall derive an HF estimate for the Riemann-Liouville integral of arbitrary order $\alpha$ ($\alpha$ can be an integer or a non-integer (real number)) of $f(t)$ via the generalized one-shot operational matrices.

Let us recall the definition of the Riemann-Liouville fractional integral of order $\alpha$ of $f(t)$.

$$J^\alpha f(t) = \frac{1}{\Gamma(\alpha)} \int_0^t (t-\tau)^{\alpha-1} f(\tau) d\tau. \tag{8}$$

Equation (8) can be written as,

$$J^\alpha f(t) = \frac{1}{\Gamma(\alpha)} \int_0^t (t-\tau)^{\alpha-1} f(\tau) d\tau = \frac{1}{\Gamma(\alpha)} \int_0^t \tau^{\alpha-1} f(t-\tau) d\tau = \frac{1}{\Gamma(\alpha)} \left(t^{\alpha-1} * f(t)\right). \tag{9}$$

Substituting the HF estimate of $f(t)$,

$$J^\alpha f(t) \approx \frac{1}{\Gamma(\alpha)} \left(t^{\alpha-1} * \left(C_S^T S_{(m)}(t) + C_T^T T_{(m)}(t)\right)\right),$$

$$= C_S^T \left(\frac{1}{\Gamma(\alpha)} \int_0^t (t-\tau)^{\alpha-1} S_{(m)}(\tau) d\tau\right) + C_T^T \left(\frac{1}{\Gamma(\alpha)} \int_0^t (t-\tau)^{\alpha-1} T_{(m)}(\tau) d\tau\right), \tag{10}$$

$$= C_S^T J^\alpha S_{(m)}(t) + C_T^T J^\alpha T_{(m)}(t),$$

$$= C_S^T \left[J^\alpha S_0(t) \ \ J^\alpha S_1(t) \ \ \cdots \ \ J^\alpha S_{m-1}(t)\right]^T + C_T^T \left[J^\alpha T_0(t) \ \ J^\alpha T_1(t) \ \ \cdots \ \ J^\alpha T_{m-1}(t)\right]^T,$$

where $[\cdots]^T$ signifies transpose.

Utilization of the orthogonal HFs for approximating $f(t)$ in (10) transforms the fractional integration of $f(t)$ to the fractional integration of the SHF set and the TF set.

Performing fractional integration on the first member of the SHF set $S_{(m)}(t)$,

$$J^\alpha S_0(t) = \frac{1}{\Gamma(\alpha)}\int_0^t (t-\tau)^{\alpha-1} S_0(\tau)d\tau = \begin{cases} 0, & \text{for } t = 0, \\ \frac{h^\alpha}{\Gamma(\alpha+1)}\left(j^\alpha - (j-1)^\alpha\right), & \text{for } t \in (0,T] \end{cases} \quad (11)$$

Evaluating (11) at $j = 1,2,3,\cdots\cdots$ provides the samples of $J^\alpha S_0(t)$ at $h$, $2h$, $3h$, $4h$, $\cdots\cdots$, $mh$. Using these samples, we can approximate $J^\alpha S_0(t)$ in the orthogonal HF domain as given in the following equation.

$$J^\alpha S_0(t) = \frac{h^\alpha}{\Gamma(\alpha+1)}\begin{bmatrix} 0 & \varsigma_1 & \varsigma_2 & \varsigma_3 & \cdots & \varsigma_{m-1} \end{bmatrix} S_{(m)}(t) + \frac{h^\alpha}{\Gamma(\alpha+1)}\begin{bmatrix} 1 & \xi_1 & \xi_2 & \xi_3 & \cdots & \xi_{m-1} \end{bmatrix} T_{(m)}(t),$$

(12)

where $\varsigma_k = \left(k^\alpha - (k-1)^\alpha\right)$, $\xi_k = (k+1)^\alpha - 2k^\alpha + (k-1)^\alpha$.

Similarly,

$$J^\alpha S_1(t) = \frac{h^\alpha}{\Gamma(\alpha+1)}\begin{bmatrix} 0 & 0 & \varsigma_1 & \varsigma_2 & \cdots & \varsigma_{m-2} \end{bmatrix} S_{(m)}(t) + \frac{h^\alpha}{\Gamma(\alpha+1)}\begin{bmatrix} 0 & 1 & \xi_1 & \xi_2 & \cdots & \xi_{m-2} \end{bmatrix} T_{(m)}(t).$$

(13)

$$\vdots$$

$$J^\alpha S_{m-2}(t) = \frac{h^\alpha}{\Gamma(\alpha+1)}\begin{bmatrix} 0 & 0 & \cdots & \cdots & 0 & \varsigma_1 \end{bmatrix} S_{(m)}(t) + \frac{h^\alpha}{\Gamma(\alpha+1)}\begin{bmatrix} 0 & \cdots & \cdots & 0 & 1 & \xi_1 \end{bmatrix} T_{(m)}(t). \quad (14)$$

$$J^\alpha S_{m-1}(t) = \frac{h^\alpha}{\Gamma(\alpha+1)}\begin{bmatrix} 0 & 0 & \cdots & \cdots & 0 & 0 \end{bmatrix} S_{(m)}(t) + \frac{h^\alpha}{\Gamma(\alpha+1)}\begin{bmatrix} 0 & \cdots & \cdots & 0 & 0 & 1 \end{bmatrix} T_{(m)}(t). \quad (15)$$

From Equations (12) to (15),

$$J^\alpha S_{(m)}(t) = P_{\alpha ss(m)} S_{(m)}(t) + P_{\alpha st(m)} T_{(m)}(t), \quad (16)$$

where

$$P_{\alpha ss(m)} = \frac{h^\alpha}{\Gamma(\alpha+1)}[\![0 \quad \varsigma_1 \quad \varsigma_2 \quad \varsigma_3 \quad \cdots \quad \varsigma_{m-1}]\!], \quad P_{\alpha st(m)} = \frac{h^\alpha}{\Gamma(\alpha+1)}[\![1 \quad \xi_1 \quad \xi_2 \quad \xi_3 \quad \cdots \quad \xi_{m-1}]\!].$$

Following the same procedure, we can estimate $J^\alpha T_{(m)}(t)$ in the orthogonal HF domain as bestowed in the next equation.

$$J^\alpha T_{(m)}(t) = P_{\alpha ts(m)}S_{(m)}(t) + P_{\alpha tt(m)}T_{(m)}(t), \tag{17}$$

where

$$P_{\alpha ts(m)} = \frac{h^\alpha}{\Gamma(\alpha+2)}[\![0 \quad \phi_1 \quad \phi_2 \quad \phi_3 \quad \cdots \quad \phi_{m-1}]\!], \quad P_{\alpha tt(m)} = \frac{h^\alpha}{\Gamma(\alpha+2)}[\![1 \quad \psi_1 \quad \psi_2 \quad \psi_3 \quad \cdots \quad \psi_{m-1}]\!],$$

$$\phi_k = k^{\alpha+1} - (k-1)^\alpha(k+\alpha), \quad \psi_k = (k+1)^{\alpha+1} - (k+1+\alpha)k^\alpha - k^{\alpha+1} + (k+\alpha)(k-1)^\alpha.$$

The HF estimate of the Riemann-Liouville fractional integral of order $\alpha$ of the function $f(t)$ is

$$\begin{aligned}\frac{1}{\Gamma(\alpha)}\int_0^t (t-\tau)^{\alpha-1}f(\tau)d\tau &\approx C_S^T\left(P_{\alpha ss(m)}S_{(m)}(t) + P_{\alpha st(m)}T_{(m)}(t)\right) + C_T^T\left(P_{\alpha ts(m)}S_{(m)}(t) + P_{\alpha tt(m)}T_{(m)}(t)\right), \\ &= \left(C_S^T P_{\alpha ss(m)} + C_T^T P_{\alpha ts(m)}\right)S_{(m)}(t) + \left(C_S^T P_{\alpha st(m)} + C_T^T P_{\alpha tt(m)}\right)T_{(m)}(t),\end{aligned} \tag{18}$$

where $P_{\alpha ss(m)}$, $P_{\alpha st(m)}$, $P_{\alpha ts(m)}$, and $P_{\alpha tt(m)}$ are the generalized one-shot operational matrices and acting as a fractional (generalized) integrator in the orthogonal HF domain.

If the order of fractional integral is one, then $P_{\alpha ss(m)} = P_{1ss(m)}$, $P_{\alpha st(m)} = P_{1st(m)}$, $P_{\alpha ts(m)} = P_{1ts(m)}$, $P_{\alpha tt(m)} = P_{1tt(m)}$ i.e. the generalized one-shot operational matrices become the classical one-shot operational matrices like the Riemann-Liouville fractional integral of $f(t)$ reduces to the first order integral of $f(t)$ when $\alpha = 1$.

We now verify the derived HF estimate of the Riemann-Liouville fractional integral of order $\alpha$ of $f(t)$. Let us take a time function $f(t) = t, t \in [0,1]$. We choose the step size as 0.125 ($m=8$).

The exact Riemann-Liouville fractional integral of $f(t)$ is

$$J^\alpha f(t) = \frac{\Gamma(2)t^{1+\alpha}}{\Gamma(2+\alpha)}. \tag{19}$$

We use $\tilde{J}^\alpha f(t)$ to denote the HF estimate of $J^\alpha f(t)$. As we notice in Table 1, the expression in (18) formulated using the generalized one-shot HF operational matrices is able to give highly accurate approximation to $J^\alpha f(t)$ even for small value of $m$. It is, therefore, confirmed that the derived generalized one-shot HF operational matrices are precise and indeed acting as a generalized integrator in the orthogonal HF domain.

It is found (Table 2) that the one-shot operational matrices presented in (Deb et al. (2012)) for $n$ ( $n$ is an integer) times repeated integration become incorrect when $n \geq 3$. So they fail to act as an integrator in the HF domain for $n \geq 3$. Using these one-shot operational matrices for solving integral or differential equations involving integrals and/or derivatives of order greater than or equal to three introduces larger error than the error produced via our approach (Equation (18)). To prove this fact, let us consider the numerical example from Section 9 of (Deb et al. (2012)).

$$f(t) = \int t\,dt + \int\int t\,dt + \int\int\int t\,dt, \ t \in [0,1], \ h = 0.125. \tag{20}$$

The approximate solution of (37) is computed by using (18) and the percentage error is calculated and compared (Table 3) to that obtained in (Deb et al. (2012)). As the third order one-shot operational matrices are inaccurate, the approximate solution in (Deb et al. (2012)) is less accurate than the approximate solution acquired by our approach.

## 4. Numerical method to solve system of fractional order differential equations

We develop here a numerical method to find the approximate solution of the system of fractional order differential equations via the orthogonal HFs.

Let us consider the following general form of the system of fractional order differential equations,

$$^C_0 D_t^\alpha y_i(t) = f_i(t, y_1(t), y_2(t), \ldots, y_n(t)), \ i = 1,2,\ldots,n, \ t \in [0,1], \ \alpha \in (0,1], \tag{21}$$

with the initial conditions are $y_i(0) = c_i$, $i = 1,2,\ldots,n$.

We describe the $i^{th}$ unknown function, $y_i(t)$, of the system of fractional order differential equations in (21) as

$$y_i(t) = y_{i,0} + z_i(t), \ i = 1,2,\ldots,n, \tag{22}$$

where $z_i(t)$ is a new unknown function, $y_{i,0}$ can be chosen such that $y_i(0) - y_{i,0} = 0$.

Employing (22) and carrying out fractional integration on (21) results in the form given below.

$$z_i(t) = J^\alpha \left( f_i(t, z_1(t) + y_{1,0}, z_2(t) + y_{2,0}, \ldots, z_n(t) + y_{n,0}) \right). \tag{23}$$

This modified system of fractional Volterra integral equations owns zero initial conditions.

Expanding the unknown function, $z_i(t)$, and the nonlinear function, $f_i(t, z_1(t) + y_{1,0}, z_2(t) + y_{2,0}, \ldots, z_n(t) + y_{n,0})$, into orthogonal HFs,

$$z_i(t) \approx C_{Si}^T S_{(m)}(t) + C_{Ti}^T T_{(m)}(t), \ \forall i \in [1,n], \tag{24}$$

where $C_{Si}^T = [c_{i0} \ c_{i1} \ c_{i2} \ \cdots \ \cdots \ c_{i(m-1)}]$, $c_{ij} = z_i(jh)$, $j = 0,1,2,\ldots,m$.

$C_{Ti}^T = [(c_{i1} - c_{i0}) \ (c_{i2} - c_{i1}) \ \cdots \ \cdots \ (c_{im} - c_{i(m-1)})]$.

$$f_i\left(t, z_1(t)+ y_{1,0}, z_2(t)+ y_{2,0},\ldots\ldots, z_n(t)+ y_{n,0}\right) \approx \tilde{C}_{Si}^T S_{(m)}(t) + \tilde{C}_{Ti}^T T_{(m)}(t), \quad \forall i \in [1,n], \quad (25)$$

where $\tilde{C}_{Si}^T = \begin{bmatrix} e_{i0} & e_{i1} & e_{i2} & \cdots & \cdots & e_{i(m-1)} \end{bmatrix}$, $\tilde{C}_{Ti}^T = \begin{bmatrix} (e_{i1}-e_{i0}) & (e_{i2}-e_{i1}) & \cdots & \cdots & (e_{im}-e_{i(m-1)}) \end{bmatrix}$

$e_{ij} = f_i\left(t, z_1(jh)+ y_{1,0}, z_2(jh)+ y_{2,0},\ldots\ldots, z_n(jh)+ y_{n,0}\right), \quad j = 0,1,2,\ldots\ldots, m$.

From (23) to (25),

$$C_{Si}^T S_{(m)}(t) + C_{Ti}^T T_{(m)}(t) = J^\alpha \left(\tilde{C}_{Si}^T S_{(m)}(t) + \tilde{C}_{Ti}^T T_{(m)}(t)\right). \quad (26)$$

Employing the generalized one-shot operational matrices,

$$C_{Si}^T S_{(m)}(t) + C_{Ti}^T T_{(m)}(t) = \tilde{C}_{Si}^T \left(P_{\alpha ss(m)} S_{(m)}(t) + P_{\alpha st(m)} T_{(m)}(t)\right) + \tilde{C}_{Ti}^T \left(P_{\alpha ts(m)} S_{(m)}(t) + P_{\alpha tt(m)} T_{(m)}(t)\right). \quad (27)$$

Equating the coefficients of SHF set, $S_{(m)}(t)$, and TF set, $T_{(m)}(t)$,

$$C_{Si}^T = \tilde{C}_{Si}^T P_{\alpha ss(m)} + \tilde{C}_{Ti}^T P_{\alpha ts(m)}, \quad C_{Ti}^T = \tilde{C}_{Si}^T P_{\alpha st(m)} + \tilde{C}_{Ti}^T P_{\alpha tt(m)}, \quad i = 1,2,\ldots\ldots, n \quad (28)$$

Solving the above system of nonlinear algebraic equations gives the coefficient vectors; $C_{S1}^T, C_{S2}^T, \ldots, C_{Sn}^T$ and $C_{T1}^T, C_{T2}^T, \ldots, C_{Tn}^T$.

From Equation (24), we have the HF estimate of the new unknown function $z_i(t)$.

The actual unknown, $y_i(t)$, is approximated in the orthogonal HF domain as

$$y_i(t) = y_{i,0} + C_{Si}^T S_{(m)}(t) + C_{Ti}^T T_{(m)}(t), \quad i \in [1,n]. \quad (29)$$

## 5. Convergence analysis

In this section, we shall prove that the HF approximate solution of the system of fractional order differential equations in (21) converges to its exact solution when $m$ is large enough.

Let $Y(t) = [y_1(t) \quad y_2(t) \quad y_3(t) \quad \cdots \quad y_n(t)]$ and $\tilde{Y}(t) = [\tilde{y}_1(t) \quad \tilde{y}_2(t) \quad \tilde{y}_3(t) \quad \cdots \quad \tilde{y}_n(t)]$ be the exact and the approximate solution of (21), respectively.

Let us define an error between $y_i(t)$ and $\tilde{y}_i(t)$ as

$$\varepsilon_i = \|y_i(t) - \tilde{y}_i(t)\|, \quad i = 1,2,\ldots\ldots, n. \quad (30)$$

Let us suppose that the nonlinear function $f_i(t, Y(t))$ satisfies the Lipschitz condition uniformly in $t \in [0,1]$. Therefore, there exists a positive constant, $L$, called Lipschitz constant such that

$$\|f_i(t, Y(t)) - f_i(t, \tilde{Y}(t))\| \leq L \|Y(t) - \tilde{Y}(t)\|, \quad i = 1,2,\ldots\ldots, n. \quad (31)$$

Under these assumptions, we state the following theorem.

**Theorem 5.1** The HF approximate solution of (21) converges to its exact solution if and only if $0 < nLT^\alpha / \Gamma(\alpha+1) < 1$.

**Proof**

Equation (30) can be written as

$$\varepsilon_i = \left\| \frac{1}{\Gamma(\alpha)} \int_0^T (t-\tau)^{\alpha-1} f_i(\tau, Y(\tau)) d\tau - \frac{1}{\Gamma(\alpha)} \int_0^T (t-\tau)^{\alpha-1} f_i(\tau, \tilde{Y}(\tau)) d\tau \right\|, \tag{32}$$

$$= \left| \frac{1}{\Gamma(\alpha)} \int_0^T (t-\tau)^{\alpha-1} \left\| f_i(\tau, Y(\tau)) - f_i(\tau, \tilde{Y}(\tau)) \right\| d\tau \right|, \tag{33}$$

$$\leq \left| \frac{1}{\Gamma(\alpha)} \int_0^T (t-\tau)^{\alpha-1} L \| Y(\tau) - \tilde{Y}(\tau) \| d\tau \right|, \tag{34}$$

$$\leq \left| \frac{L}{\Gamma(\alpha)} \int_0^T (t-\tau)^{\alpha-1} \sum_{i=1}^n |y_i(\tau) - \tilde{y}_i(\tau)| d\tau \right|, \tag{35}$$

$$\leq \left| \frac{Ln\varepsilon_i}{\Gamma(\alpha)} \int_0^T (t-\tau)^{\alpha-1} d\tau \right|, \tag{36}$$

$$\leq \frac{T^\alpha Ln\varepsilon_i}{\Gamma(\alpha+1)}. \tag{36}$$

Therefore,

$$\varepsilon_i \left( 1 - \frac{T^\alpha Ln}{\Gamma(\alpha+1)} \right) \leq 0 \tag{37}$$

If $\left( T^\alpha nL / \Gamma(\alpha+1) \right) \in (0,1)$, the error $\varepsilon_i$ tends to zero when sufficiently large number of subintervals is considered. This completes the proof.

## 6. Numerical examples

In this section, the numerical method devised in Section 4 is applied to the systems of linear and nonlinear fractional order differential equations. The Matlab built-in function 'fsolve' is used to solve the system of nonlinear algebraic equations arisen upon the utilization of the generalized one-shot operational matrices in place of the Riemann-Liouville fractional order integral.

### *Example 6.1*

Consider the following system of nonlinear fractional order differential equations.

$$D^\alpha y_1(t) = y_1(t) + y_2^2(t), \quad y_1(0) = 0, \quad y_1'(0) = 1, \tag{38}$$

$$D^\beta y_2(t) = y_1(t) + 5y_2(t), \ y_2(0)=0, \ y_2'(0)=1, \ y_2''(0)=1. \tag{39}$$

The analytical solution is unknown.

The obtained piecewise linear HF solutions (Subplot (a) in Figure 1) by our numerical method using the step size of 0.001 for $\alpha = 1.3$ and $\beta = 2.4$ are in good accordance with the solutions acquired by Adomian decomposition method in [Varsha and Hossein 2007], via fractional differential transform method in [Ertürk and Momani 2008] (see Example 2 in [Ertürk and Momani 2008]) and by Legendre wavelets based numerical method in [Yiming et al. 2015] (see Example 5.2 in [Yiming et al. 2015]). For different $\alpha$ and $\beta$, our numerical method (with $h = 0.001$) produced the same results (Subplots (b) and (c) in Figure 1) as Legendre wavelets based numerical method provided in [Yiming et al. 2015] (see Example 5.2 in [Yiming et al. 2015]).

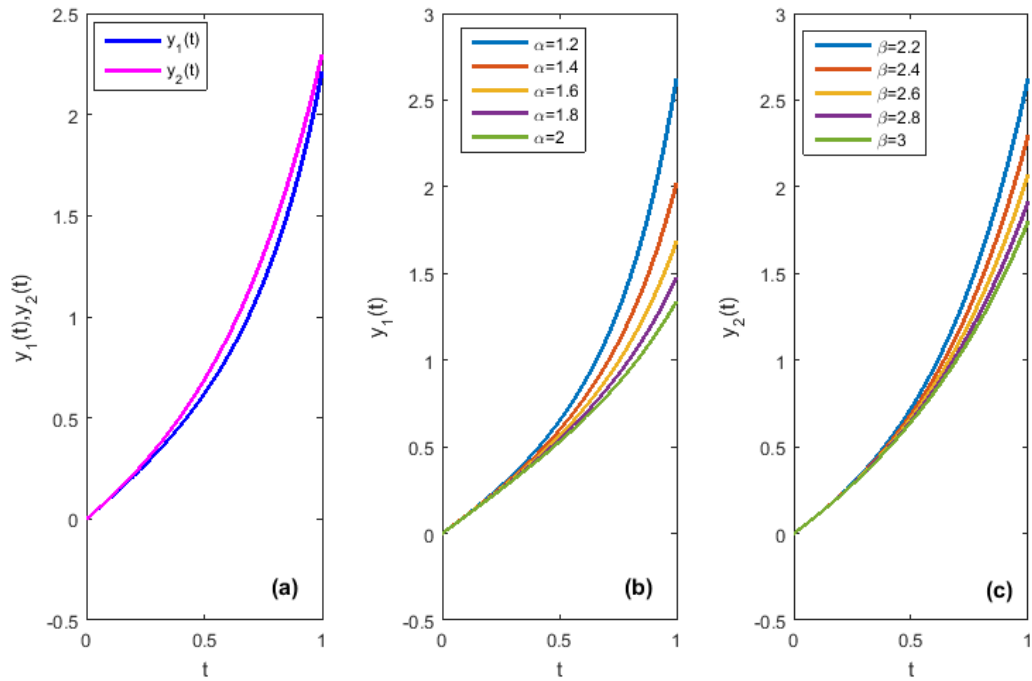

**Figure 1** HF solutions of Example 6.1 for $\alpha = 1.3$ and $\beta = 2.4$ (Subplot (a)) and for other values of $\alpha$ and $\beta$ (Subplots (b) and (c))

*Example 6.2*

Consider the system of linear fractional order differential equation.

$$D^\alpha x(t) = x(t) + y(t), \ x(0) = 0. \tag{40}$$

$$D^\beta y(t) = -x(t) + y(t), \ y(0) = 1. \tag{41}$$

The exact solution when $\alpha = 1$ and $\beta = 1$ is $x(t) = e^t \sin t$, $y(t) = e^t \cos t$.

The absolute error between the exact solution and the HF solution is computed for different values of $h$ and displayed in Figure 2. Table 4 presents the $\infty$-norm of the absolute error. The HF solutions for the integer order case ($\alpha = 1$, $\beta = 1$) and the fractional order case ($\alpha = 0.7$, $\beta = 0.9$) shown in Figure 3 are in good agreement with the approximate solutions obtained in [Ertürk and Momani 2008] (see Example 1 in [Ertürk and Momani 2008]), [Khader et al. 2015] (see Example 1 in [Khader et al. 2015]), [Momani and Odibat 2007] (see Example 6.2 in [Momani and Odibat 2007]) and [Zurigat et al. 2010] (see Example 4.1 in [Zurigat et al. 2010]).

**Table 4** Error analysis of example 6.2

| $h$ | $\|e_1\|_\infty$ | $\|e_2\|_\infty$ |
| --- | --- | --- |
| 1/10 | 0.001387236644377 | 0.006249545001395 |
| 1/200 | 3.411247368134700e-06 | 1.565014461890610e-05 |
| 1/400 | 8.527964938664920e-07 | 3.912552036577920e-06 |
| 1/600 | 3.790271057013680e-07 | 1.738915393678650e-06 |
| 1/800 | 2.132093102069630e-07 | 9.781421581589460e-07 |
| 1/1000 | 1.364585999752420e-07 | 6.260121305778910e-07 |

*Example 6.3*

The system of nonlinear fractional order differential equations is

$$D^\alpha x(t) = 2y^2(t), \quad x(0) = 0. \tag{42}$$

$$D^\beta y(t) = tx(t), \quad y(0) = 1. \tag{43}$$

$$D^\gamma z(t) = y(t)z(t), \quad z(0) = 1. \tag{44}$$

The given problem has no closed form solution.

For $\alpha = \beta = \gamma = 1$ (Subplot (a) in Figure 4) and $\alpha = 0.8, \beta = 0.7, \gamma = 0.6$ (Subplot (b) in Figure 4), our numerical method (with step size of 0.001) produced the same results as attained in [Khader et al. 2015] (see Example 2 in [Khader et al. 2015]) and [Chen et al. 2015] (see Example 5.3 in [Chen et al. 2015]). The piecewise linear HF solution (Subplot (c) in Figure 4) for $\alpha = 0.75, \beta = 0.85, \gamma = 0.95$ is in accordance with the solution of Adomian decomposition method and variational iteration method in [Momani and Odibat 2007] (see Example 6.4 in [Momani and Odibat 2007]). Figure 5 shows the results of our method for other values of $\alpha$, $\beta$ and $\gamma$, which exactly match the solutions obtained in [Chen et al. 2015].

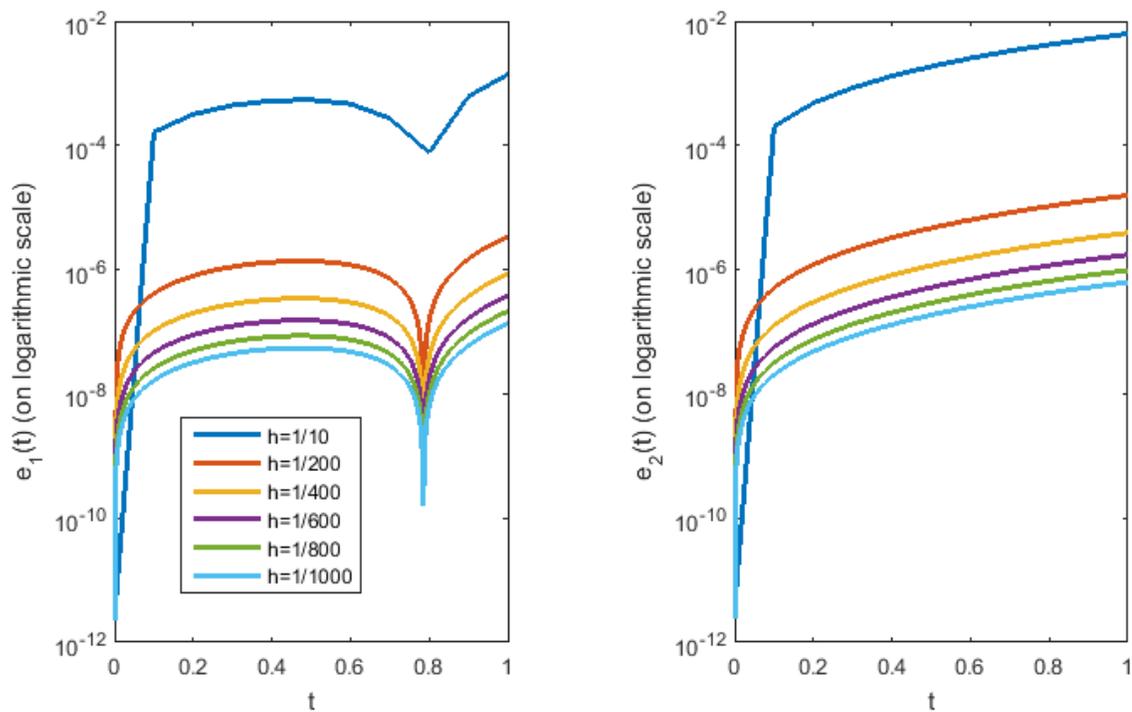

**Figure 2** Absolute error produced via our numerical method for $\alpha = 1$

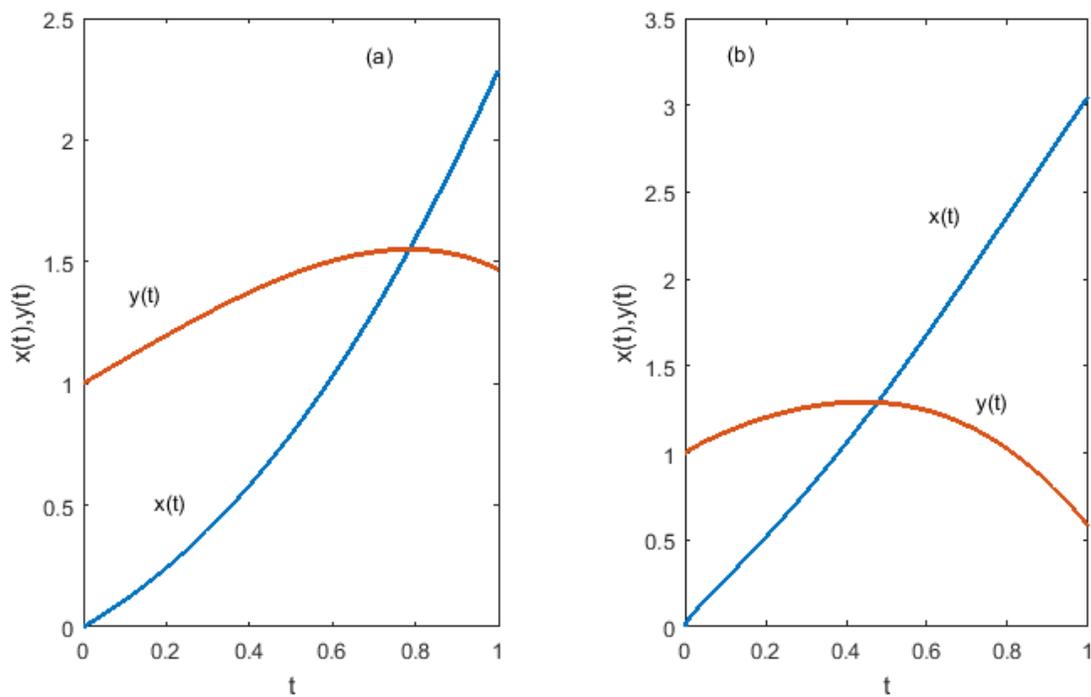

**Figure 3** HF solution of Example 6.2 for $\alpha = \beta = 1$ (a) and $\alpha = 0.7, \beta = 0.9$ (b)

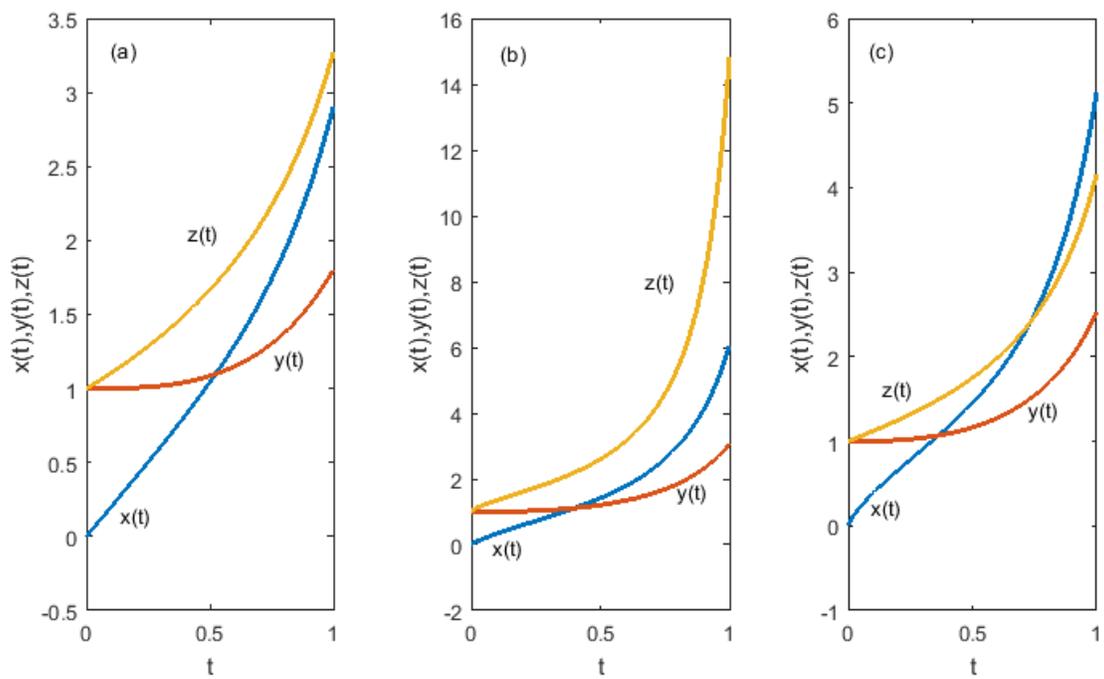

**Figure 4** HF solution of Example 6.3 for $\alpha = \beta = \gamma = 1$ (a), $\alpha = 0.8, \beta = 0.7, \gamma = 0.6$ (b) and $\alpha = 0.75, \beta = 0.85, \gamma = 0.95$ (c)

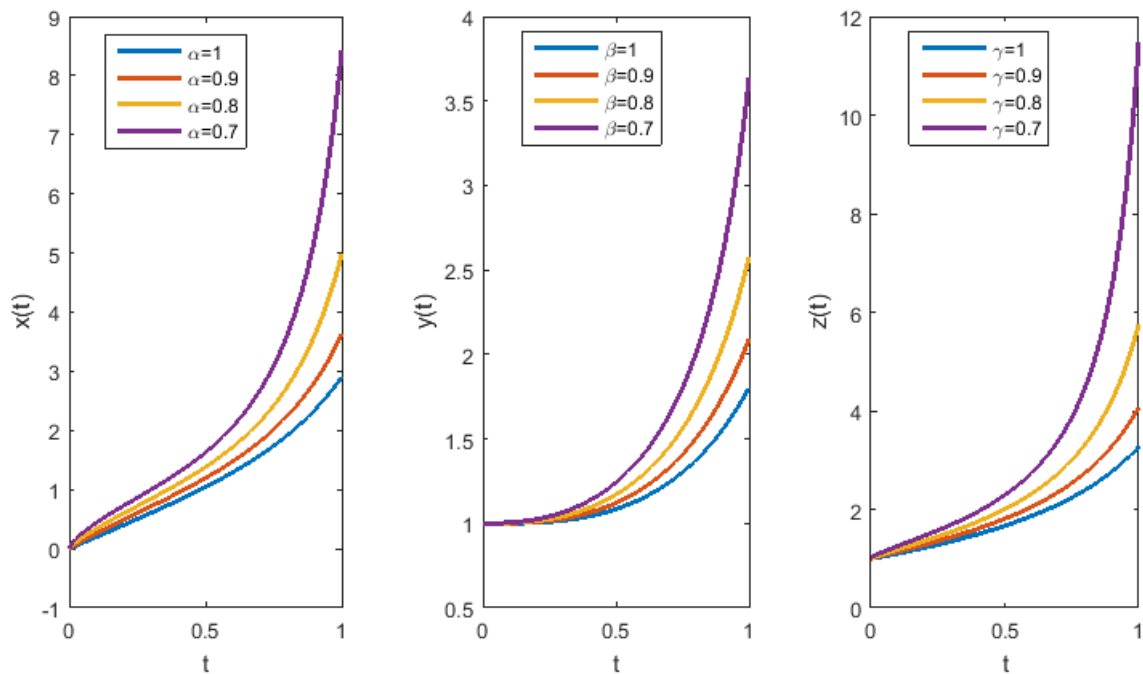

**Figure 5** HF solutions of Example 6.3 for other values of $\alpha$, $\beta$ and $\gamma$

**Example 6.4** (Muhaya (2013))

The fractional order mathematical model describing the effects of smoking in a population is given as

$$^C_0D^\alpha_t S_p(t) = (1-p)\pi - (\varepsilon + \lambda + \mu)S_p(t), \tag{45}$$

$$^C_0D^\alpha_t S_a(t) = p\pi + \varepsilon S_p(t) - (1-f)\lambda S_a(t) - \mu S_a(t), \tag{46}$$

$$^C_0D^\alpha_t L_S(t) = \lambda S_p(t) + (1-f)\lambda S_a(t) + \sigma C_S(t) - (\mu + \rho_1 + \gamma)L_S(t), \tag{47}$$

$$^C_0D^\alpha_t C_S(t) = \gamma L_S(t) - (\sigma + \rho_2 + \delta + \mu)C_S(t), \tag{48}$$

$$^C_0D^\alpha_t Q(t) = \rho_1 L_S(t) + \rho_2 C_S(t) - \mu Q(t), \tag{49}$$

where $\lambda = \beta(L_S(t) + \eta C_S(t))/N(t)$, $N(t) = S_p(t) + S_a(t) + L_S(t) + C_S(t) + Q(t)$, $t \in [0,1]$, $\alpha \in (0,1]$.

Here $S_p(t)$ is the number of potential smokers reacting to pro-smoking advertisements, $S_a(t)$ is the number of potential smokers responding to anti-smoking ads, $L_S(t)$ is the number of light smokers, $C_S(t)$ is the number of chain smokers and $Q(t)$ is the number of quitters.

We use the following model parameters and initial conditions.

$\rho_1 = 0.5$, $\rho_2 = 0.25$, $\varepsilon = 0.001$, $\sigma = 0.0307$, $p = 0.8$, $f = 1$, $\beta = 2$, $\pi = 14$, $\mu = 0.031$, $\delta = 0.01$, $\eta = 0.0002$, $\gamma = 0.6$, $S_p(0) = 8000$, $S_a(0) = 1970$, $L_S(0) = 20$, $C_S(0) = 10$, $Q(0) = 0$.

The fractional order smoking model is a bit complex than the fractional order model bestowed in the preceding subsection. Figure 6 displays the comparison of HF solution and the solution via RK fourth order method and Figure 7 presents the numerical solutions of fractional order smoking model obtained for different value of $\alpha$. The step size of 0.002 is employed for performing all numerical simulations. Despite the complexity and nonlinearity, the proposed numerical method demonstrates its competence to solve high dimensional system of nonlinear fractional order differential equations.

**Example 6.5** (Acevedo-Estefania et al. (2000))

The fractional order mathematical model for lung cancer is

$$^C_0D^\alpha_t N(t) = (1-q)\Lambda - \frac{\beta N(t)(I_1(t) + I_2(t))}{T(t)} - \mu N(t), \tag{50}$$

$$^C_0D^\alpha_t I_1(t) = \frac{((1-p_n)\beta N(t) + (1-p_s)\beta S(t))(I_1(t) + I_2(t))}{T(t)} - (\sigma_1 + \gamma_1 + \delta_1 + \mu)I_1(t), \tag{51}$$

$$_0^C D_t^\alpha I_2(t) = \gamma_1 I_1(t) - (\gamma_2 + \delta_2 + \mu) I_2(t), \tag{52}$$

$$_0^C D_t^\alpha Q(t) = p_2 \gamma_2 I_2(t) + p_1 \sigma_1 I_1(t) - (\delta_q + \mu) Q(t), \tag{53}$$

$$_0^C D_t^\alpha S(t) = (1 - p_1) \sigma_1 I_1(t) + (1 - p_2) \gamma_2 I_2(t) - \frac{\beta S(t)(I_1(t) + I_2(t))}{T(t)} - \mu S(t), \tag{54}$$

$$_0^C D_t^\alpha L(t) = \frac{(p_n \beta N(t) + p_s \beta S(t) + \beta_e E(t))(I_1(t) + I_2(t))}{T(t)} + \delta_1 I_1(t) + \delta_2 I_2(t) + \delta_q Q(t) - (\mu + d) L(t), \tag{56}$$

$$_0^C D_t^\alpha E(t) = q\Lambda - \frac{\beta_e E(t)(I_1(t) + I_2(t))}{T(t)} - \mu E(t), \tag{57}$$

where $T(t) = N(t) + I_1(t) + I_2(t) + Q(t) + S(t) + L(t) + E(t)$, $t \in [0,1]$, $\alpha \in (0,1]$.

In the above fractional model, $N(t)$ signifies the number of non-smokers, $I_1(t)$ the number of light smokers, $I_2(t)$ the number of heavy smokers, $Q(t)$ individuals who quit smoking permanently, $S(t)$ individuals who stop smoking temporarily, $L(t)$ the number of smokers who develop lung cancer and $E(t)$ the educated individuals.

The following data is considered for numerical simulations.

$\Lambda = 14$, $\mu = 0.014$, $\beta = 2$, $p_n = 10^{-4}$, $\gamma_1 = 0.6$, $\gamma_2 = 0.25$, $p_1 = 0.025$, $p_2 = 0.025$, $\sigma_1 = 0.5$, $p_s = 10^{-3}$, $\delta_q = 0.005$, $\delta_2 = 0.03$, $d = 0.016$, $q = 0.25$, $\beta_e = 10^{-4}$, $N(0) = 500$, $I_1(0) = 200$, $\delta_1 = 0.01$, $I_2(0) = 200$, $Q(0) = 200$, $S(0) = 200$, $L(0) = 200$, $E(0) = 200$, $h = 0.002$.

Figures 8 and 9 express the fact that the proposed numerical method is so powerful that it can solve even more complicated and higher dimensional system of nonlinear fractional order differential equations.

**Example 6.6** (Abdulrahman et al. (2013))

The generalized mathematical model of Hepatitis B infection is

$$_0^C D_t^\alpha S_U(t) = bN(t)(1 - \varepsilon_\rho \tau_b) - b\theta(A_F(t) + \phi C_F(t))(1 - \varepsilon_\rho \tau_b) - (\sigma_s + \varepsilon_\rho \tau_U + \mu) S_U(t), \tag{58}$$

$$_0^C D_t^\alpha S_F(t) = \sigma_s S_U(t) + wV(t) - \frac{pc(A_F(t) + \eta C_F(t))(1 - \varepsilon_c \tau_c)}{N(t)} S_F(t) - (\varepsilon_p \tau_F + \mu) S_F(t), \tag{59}$$

$$_0^C D_t^\alpha V(t) = bN(t)\varepsilon_\rho \tau_b + \varepsilon_\rho \tau_U S_U(t) + \varepsilon_\rho \tau_F S_F(t) - (w + \mu) V(t), \tag{60}$$

$$_0^C D_t^\alpha A_U(t) = b\theta(A_F(t) + \phi C_F(t))(1 - \varepsilon_\rho \tau_b) - (\sigma_A + \mu + \delta_A) A_U(t), \tag{61}$$

$$^C_0D^\alpha_t A_F(t) = \frac{pc(A_F(t) + \eta C_F(t))(1-\varepsilon_c\tau_c)}{N(t)} S_F(t) - (\sigma_A + \mu + \delta_A)A_F(t), \tag{62}$$

$$^C_0D^\alpha_t C_U(t) = \sigma_A \varphi_U A_U(t) - (\sigma_c + \mu)C_U(t), \tag{63}$$

$$^C_0D^\alpha_t C_F(t) = \sigma_A \varphi_F A_F(t) + \sigma_C C_U(t) - (\gamma_C + \mu + \delta_C)C_F(t), \tag{64}$$

$$^C_0D^\alpha_t R(t) = \sigma_A(1-\varphi_U)A_U(t) + \sigma_A(1-\varphi_F)A_F(t) + \gamma_C C_F(t) - \mu R(t), \tag{65}$$

$$N(t) = S_U(t) + S_F(t) + V(t) + A_U(t) + A_F(t) + C_U(t) + C_F(t) + R(t),\ t \in [0,1],\ \alpha \in (0,1], \tag{66}$$

where $S_U(t)$ is the susceptible individuals under 15 years of age, $S_F(t)$ is susceptible individuals at or above 15 years of age, $V(t)$ is vaccinated individuals, $A_U(t)$ is acutely infected individuals under 15 years of age, $A_F(t)$ is acutely infected individuals at or above 15 years of age, $C_U(t)$ is chronically infected individuals under 15 years of age, $C_F(t)$ is chronically infected individuals at or above 15 years of age and $R(t)$ is removed individuals due to recovery from infection.

The model parameters and the initial conditions are

$b = 0.036$, $\mu = 0.021$, $c = 20$, $p = 0.079$, $\eta = 0.667$, $\phi = 0.159$, $\sigma_S = 0.067$, $\sigma_A = 2.667$, $\sigma_C = 0.068$, $\varphi_U = 0.885$, $\varphi_F = 0.1$, $\gamma_C = 0.015$, $\delta_A = 0.007$, $\delta_C = 0.001$, $\varepsilon_\rho = 0.9$, $w = 0.04$, $\varepsilon_C = 0.8$, $\tau_b = 0.66$, $\tau_C = 0.2$, $\tau_U = 0.001$, $\tau_F = 0.001$, $\theta = 0.724$, $S_U(0) = 23148265$, $S_F(0) = 51967535$, $V(0) = 16528817$, $A_U(0) = 3012259$, $A_F(0) = 5280601$, $C_U(0) = 5394338$, $C_F(0) = 6315313$, $R(0) = 40570172$.

The fractional order model of Hepatitis B infection is the most intricate of four problems we have solved so far in this paper. Using the step size of 0.002, the piecewise linear HF approximate solution is obtained via the proposed numerical method for $\alpha = 1$ (Figure 10) and $\alpha \in (0,1)$ (Figure 11). It is manifestly proved in Figure 10 that the proposed numerical method competes with the well-known numerical technique; Runge-Kutta fourth order method by offering acceptable numerical solution to such a high dimensional nonlinear system (Equations (58) to (66)). The proposed numerical method can, therefore, be used to solve highly nonlinear and high dimensional system of ordinary differential equations of integer order representing complex biological processes. In case of non-integer values of $\alpha$, the developed numerical method yields correct solutions that proves its suitability to real complex fractional order mathematical models.

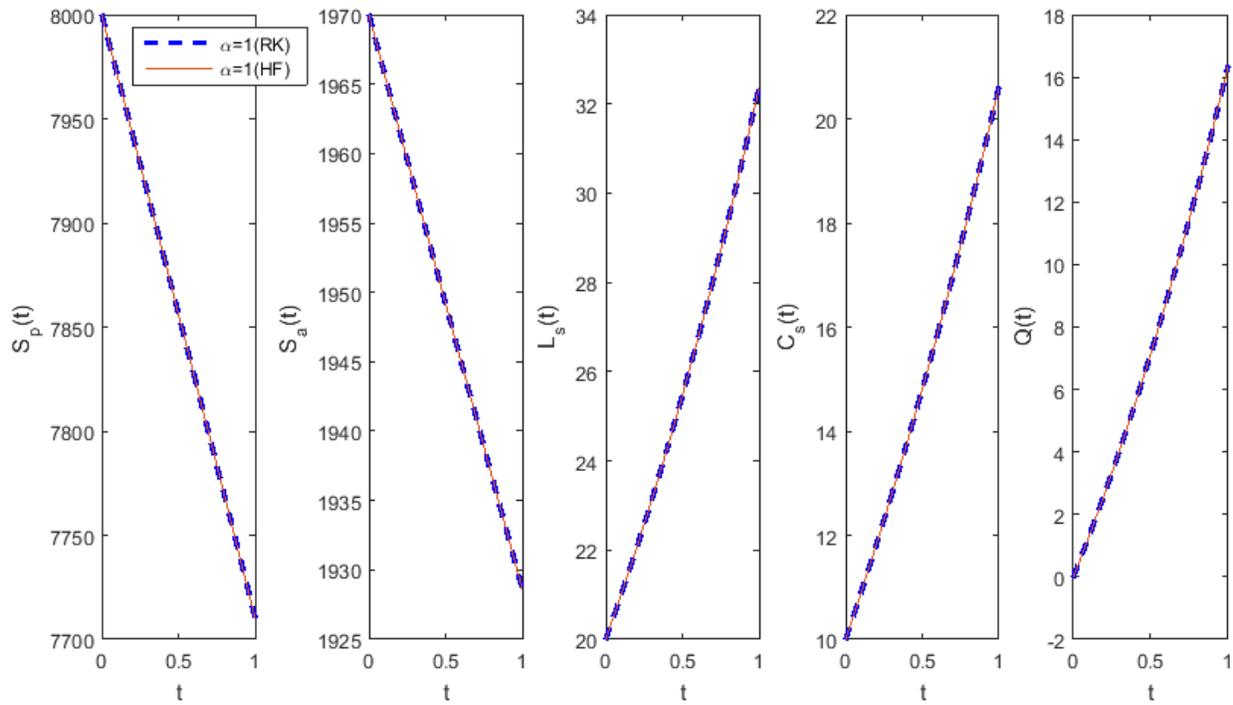

**Figure 6**: Numerical solution of smoking model of order 1.

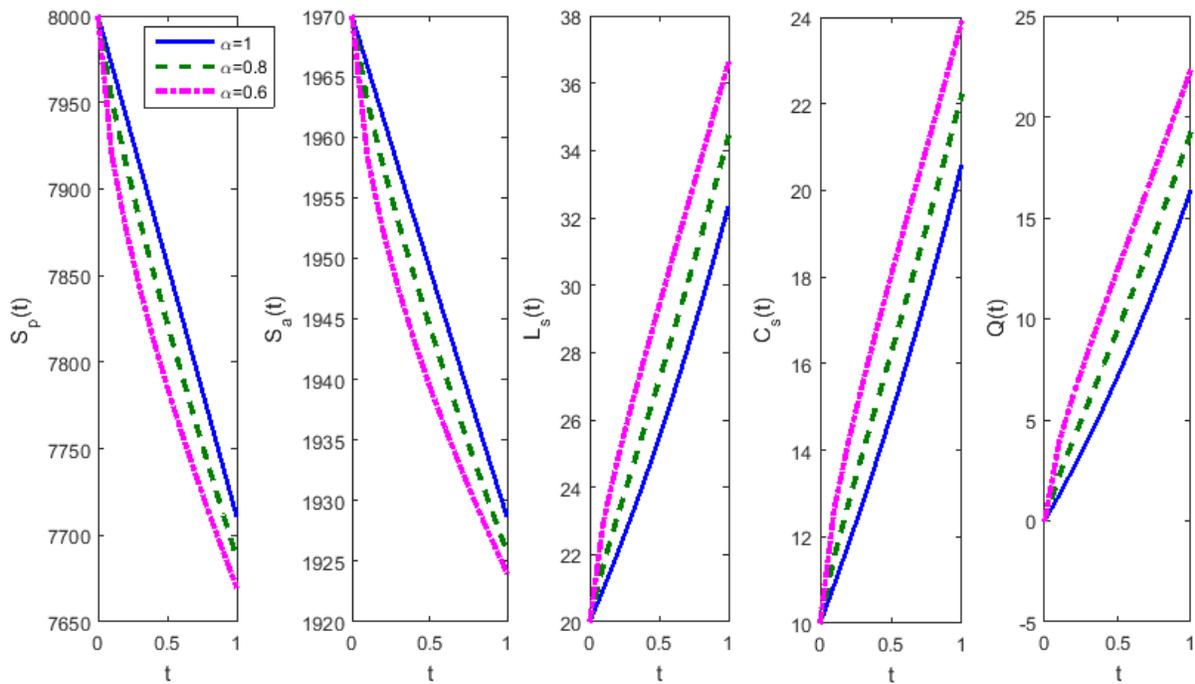

**Figure 7**: Piecewise linear HF solution of fractional order smoking model.

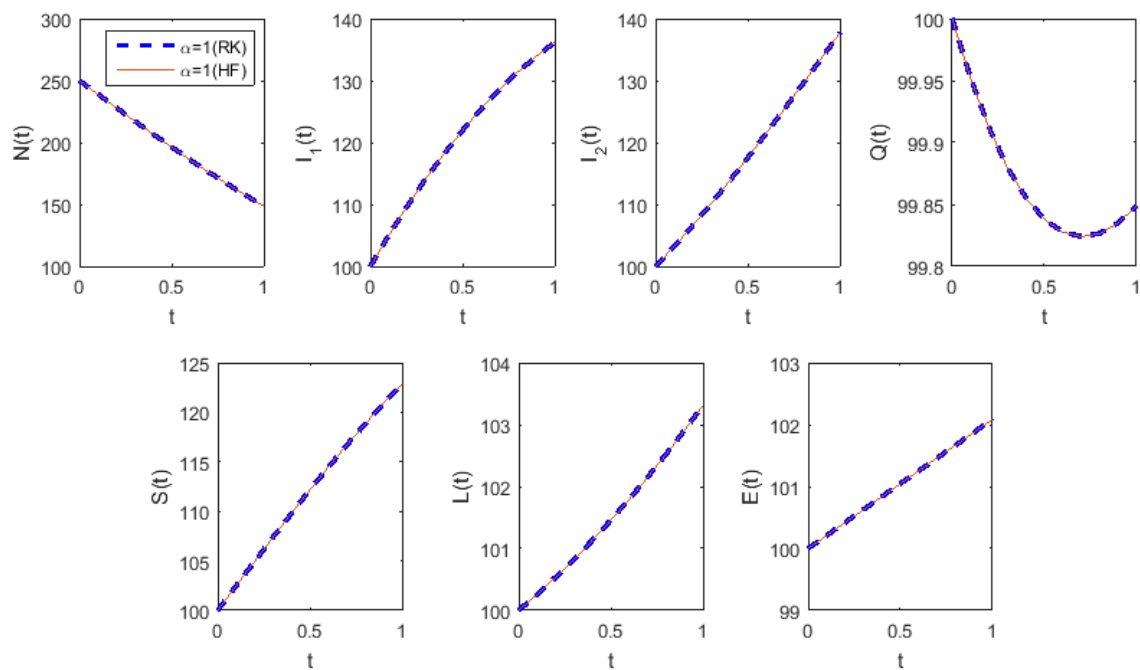

**Figure 8**: Comparison of approximate solutions of integer order ($\alpha = 1$) model for lung cancer. HF solution—solid line, RK solution—dashed line.

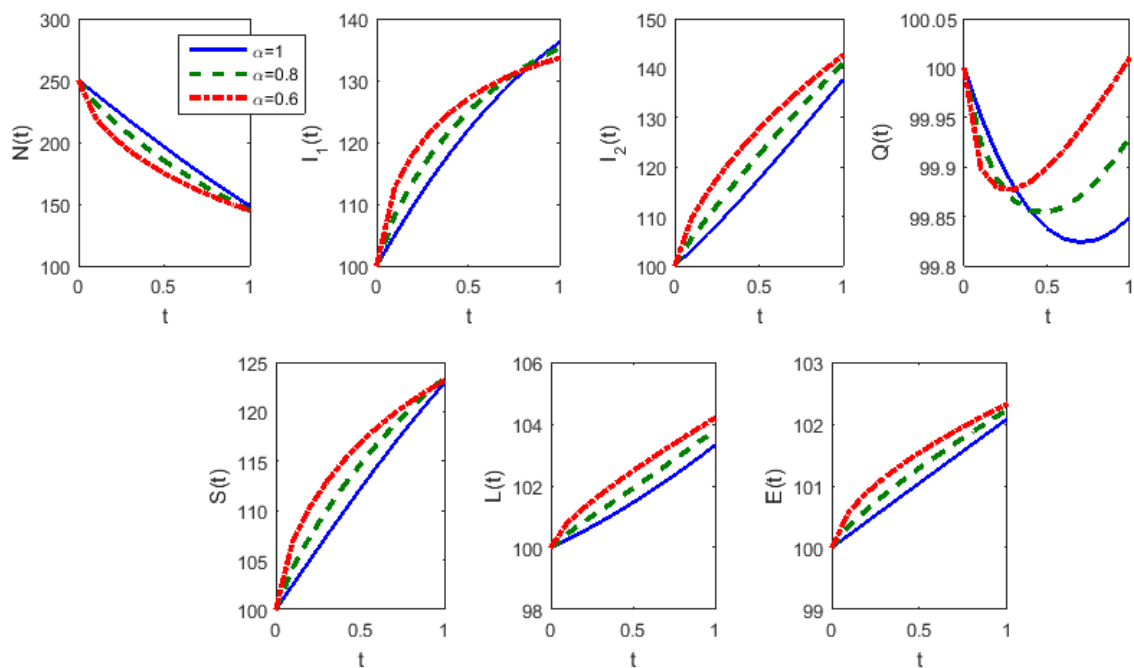

**Figure 9**: Comparison of piecewise linear HF solutions of fractional order model for lung cancer.

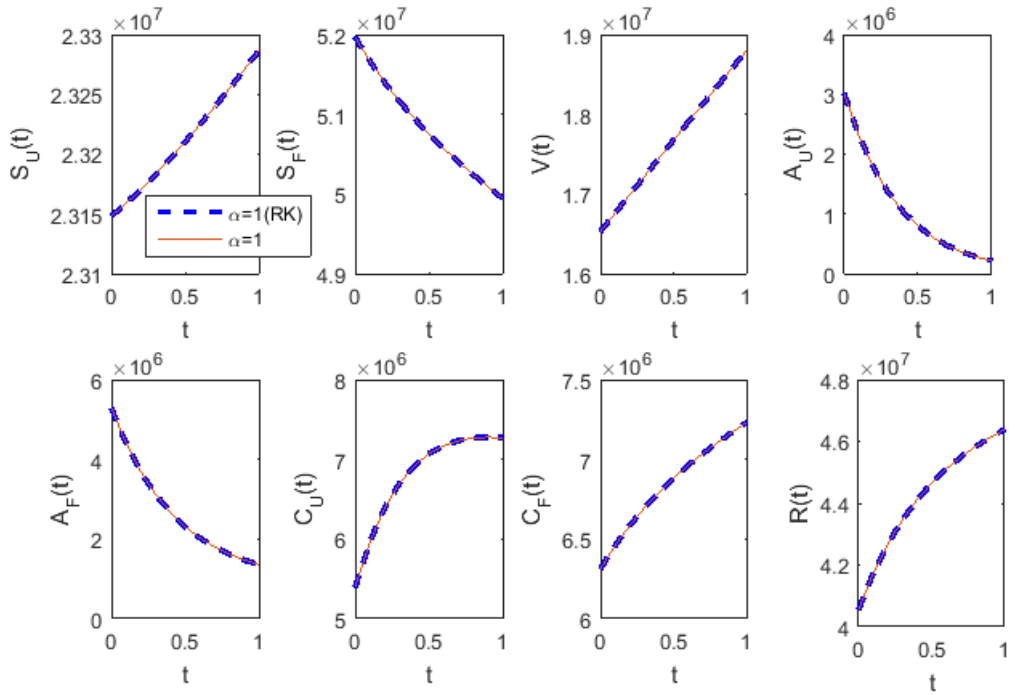

**Figure 10**: Numerical solution of integer order model ($\alpha = 1$) of Hepatitis B infection. HF solution--solid line, RK solution--dashed line.

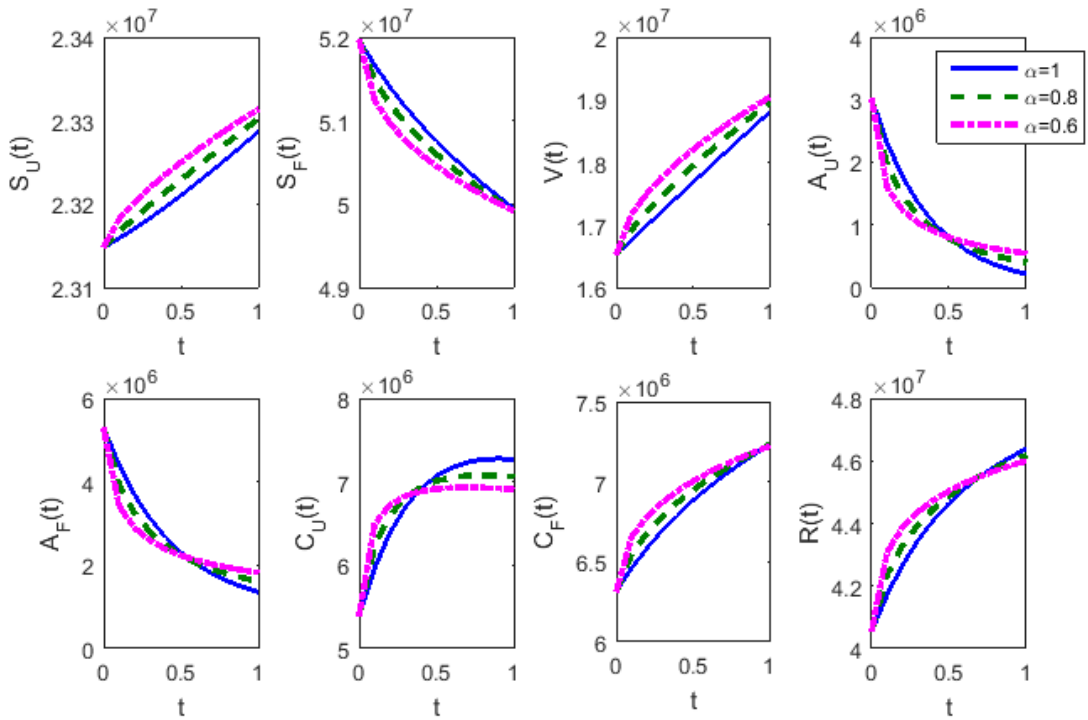

**Figure 11**: Piecewise linear HF solutions of fractional order model of Hepatitis B infection.

## 7. Conclusions

This article proposes the orthogonal hybrid function based numerical method for solving the system of differential equations of arbitrary order from the perspective of mathematics and in the pursuit of its application by applied scientist, engineers and medical practitioners. Orthogonal sample-and hold function in conjunction with triangular function was put forward as an orthogonal hybrid function and the generalized one-shot HF operational matrices derived in this work are capable of approximating integral of any order. In Examples 8.1 to 8.4, the proposed numerical method performed equally to Runge-Kutta fourth order method when $\alpha = 1$, so, it can be considered as a good alternative to Runge-Kutta fourth order method for solving integer order mathematical models numerically. Results of Section 8 and Tables 1 to 4 qualify the proposed numerical method to be a valuable tool for applied mathematicians and applied scientists to thoroughly understand innately complicated real processes. In this paper, we have successfully extended the horizon of application of the orthogonal hybrid functions.

**Appendix 1** Basic properties of HFs.

The components of SHF vector, $S_{(m)}(t)$, and TF vector, $T_{(m)}(t)$, have orthogonal property.

$$\int_0^T S_i(t)S_j(t)dt = \begin{cases} h, & \text{if } i == j, \\ 0, & \text{otherwise}, \end{cases} \int_0^T T_i(t)T_j(t)dt = \begin{cases} \dfrac{h}{3}, & \text{if } i == j, \\ \dfrac{h}{6}, & \text{otherwise}, \end{cases} \quad i \in [0, m-1]. \tag{A1}$$

The product, $S_i(t)S_j(t)$, where $i, j \in [0, m-1]$, is expressed via HFs

$$S_i(t)S_j(t) = \begin{cases} S_i(t), & \text{if } i == j, \\ 0, & \text{otherwise}. \end{cases} \tag{A2}$$

Similarly,

$$T_i(t)T_j(t) = \begin{cases} T_i(t), & \text{if } i == j, \\ 0, & \text{otherwise}, \end{cases} i, j \in [0, m-1]. \tag{A3}$$

Whereas the product, $S_i(t)T_j(t)$, is estimated in the orthogonal HF domain as

$$S_i(t)T_j(t) = \begin{cases} T_j(t), & \text{if } i == j, \\ 0, & \text{othwersie}, \end{cases} i, j \in [0, m-1]. \tag{A4}$$

The product of two functions $h(t) = f_1(t)f_2(t)$ can be approximated by HFs as given in the next equation.

$$\begin{aligned} f_1(t)f_2(t) &\approx \left(C_{S1}^T S_{(m)}(t) + C_{T1}^T T_{(m)}(t)\right)\left(C_{S2}^T S_{(m)}(t) + C_{T2}^T T_{(m)}(t)\right) \\ &= \left(C_{S1}^T .* C_{S2}^T\right) S_{(m)}(t) + \left(C_{S1}^T .* C_{T2}^T + C_{T1}^T .* C_{S2}^T + C_{T1}^T .* C_{T2}^T\right) T_{(m)}(t) \end{aligned} \tag{A5}$$

The $n^{th}$ power of function, $g(t)$, ($g(t) \in C[0,T]$) is expanded into orthogonal HFs using the following expression.

$$(g(t))^n \approx C_S^T S_{(m)}(t) + C_T^T T_{(m)}(t), \tag{A6}$$

where $C_S^T = [c_0 \quad c_1 \quad \cdots \quad \cdots \quad c_{m-1}]$, $C_T^T = [d_0 \quad d_1 \quad \cdots \quad \cdots \quad d_{m-1}]$, $c_i = (g(ih))^n$.

**Table 1**: Absolute error produced by HFs

| $t$ | $\left|J^{0.5}f(t)-\tilde{J}^{0.5}f(t)\right|$ | $\left|J^{1}f(t)-\tilde{J}^{1}f(t)\right|$ | $\left|J^{1.5}f(t)-\tilde{J}^{1.5}f(t)\right|$ | $\left|J^{2}f(t)-\tilde{J}^{2}f(t)\right|$ |
|---|---|---|---|---|
| 0 | 0 | 0 | 0 | 0 |
| 0.125 | 6.938893e-18 | 0 | 0 | 0 |
| 0.25 | 1.38777e-17 | 0 | 1.73472347e-18 | 0 |
| 0.375 | 0 | 0 | 0 | 0 |
| 0.5 | 5.55111e-17 | 0 | 6.93889390e-18 | 0 |
| 0.625 | 1.11022e-16 | 0 | 1.38777878e-17 | 0 |
| 0.75 | 1.11022e-16 | 0 | 2.77555756e-17 | 0 |
| 0.875 | 1.11022e-16 | 0 | 5.55111512e-17 | 1.3877787807e-17 |
| 1 | 1.11022e-16 | 0 | 5.551115123e-17 | 2.7755575615e-17 |

**Table 2**: Comparison of accuracy of our approach and approach in (Deb et al. (2012))

| $t$ | $\left|J^{3}f(t)-\tilde{J}^{3}f(t)\right|$ | | $\left|J^{4}f(t)-\tilde{J}^{4}f(t)\right|$ | |
|---|---|---|---|---|
| | Deb et al. (2012) | Our approach | Deb et al. (2012) | Our approach |
| 0 | 0 | 0 | 0 | 0 |
| 0.125 | 0 | 0 | 0 | 0 |
| 0.25 | 4.0690e-05 | 0 | 3.814e-06 | 0 |
| 0.375 | 0.0002441 | 0 | 2.924e-05 | 0 |
| 0.5 | 0.000732 | 4.336e-19 | 0.000109 | 0 |
| 0.625 | 0.001627 | 8.673e-19 | 0.000292 | 0 |
| 0.75 | 0.003051 | 0 | 0.000642 | 0 |
| 0.875 | 0.005126 | 0 | 0.001237 | 0 |
| 1 | 0.007975 | 6.938e-18 | 0.002171 | 0 |

**Table 3**: Comparison of % error produced by two approaches

| $t$ | % error by approach in (Deb et al. (2012)) | % error by our approach |
|---|---|---|
| 0 | 0 | 0 |
| 0.125 | 0 | 0 |
| 0.25 | 0 | 0 |
| 0.375 | 0.2503 | 0 |
| 0.5 | 0.4717 | 0 |
| 0.625 | 0.7013 | 0 |
| 0.75 | 0.8226 | 0 |
| 0.875 | 0.9828 | 0 |
| 1 | 1.1153 | 1.56e-14 |